\documentclass[reqno]{gokova}
\usepackage{geompsfi,graphicx}

\def\E{\ifmmode{\mathbb E}\else{$\mathbb E$}\fi} 
\def\N{\ifmmode{\mathbb N}\else{$\mathbb N$}\fi} 
\def\R{\ifmmode{\mathbb R}\else{$\mathbb R$}\fi} 
\def\Q{\ifmmode{\mathbb Q}\else{$\mathbb Q$}\fi} 
\def\C{\ifmmode{\mathbb C}\else{$\mathbb C$}\fi} 
\def\H{\ifmmode{\mathbb H}\else{$\mathbb H$}\fi} 
\def\Z{\ifmmode{\mathbb Z}\else{$\mathbb Z$}\fi} 
\def\P{\ifmmode{\mathbb P}\else{$\mathbb P$}\fi} 
\def\T{\ifmmode{\mathbb T}\else{$\mathbb T$}\fi} 
\def\SS{\ifmmode{\mathbb S}\else{$\mathbb S$}\fi} 
\def\DD{\ifmmode{\mathbb D}\else{$\mathbb D$}\fi} 

\renewcommand{\a}{\alpha}
\renewcommand{\b}{\beta}
\renewcommand{\d}{\delta}

\newcommand{\e}{\varepsilon}
\newcommand{\g}{\gamma}

\newcommand{\s}{\sigma}

\renewcommand{\t}{\tau}

\newtheorem{thm}{Theorem}[section]

\newtheorem{lem}[thm]{Lemma}

\theoremstyle{definition}

\theoremstyle{remark}
\newtheorem{rem}{Remark}[section]

\begin{document}

\title{A Fake  Cusp  and a Fishtail}
\author[AKBULUT]{Selman Akbulut
\\[.2cm] \small{ Dedicated to Robion Kirby on the occasion of his
60'th birthday}}

\thanks{Partially supported 
by NSF grant DMS-9626204 and TUBITAK Feza Gursey Institute}

\address{Michigan State University, 
Dept of Math, E.Lansing, MI, 48824, USA}
\email{akbulut@math.msu.edu}

\volume{5}

\begin{abstract}
We construct smooth $4$-manifolds that are homeomorphic but not 
diffeomorphic to the ``cusp" and the ``fishtail", which are 
thickened singular $2$-spheres. 
\end{abstract}

\maketitle




Even though many fake
smoothings of $4$-manifolds are known to exists, we know little about the basic
building blocks of exotic smooth manifolds. This is mainly due to
the fact that we still don't know if basic manifolds like $S^{4}$,
$S^{2}\times S^{2}$, and $S^{1}\times S^{3}$ could admit fake smooth
structures. Unable to show this, we demonstrate fake smooth structures
on manifolds that are in a way ``small deformations" of $ S^{2}\times
B^{2} $ (as in \cite{A1},\cite{A2}). 

Fishtail $F$ is the tubular neighborhood of an immersed
$2$-sphere in $S^{4}$ with one self intersection; as an handlebody it
can be described as a $4$-ball with $1 $ and $2$-handles attached as in
first picture of Figure 1. Recalling the $1$-handle notation of
\cite{A0},  $F$ is obtained by removing a tubular neighborhood of the
obvious disc
$ B^{2}\subset S^{2}\times B^{2} $ (which the ``circle with dot" bounds).
 Cusp $C$ is a $4$-ball with a 2-handle
attached along a trefoil knot with $0$-framing as shown in the second
picture of Figure 1. 

\begin{figure}[htb]
\includegraphics[scale=.5]{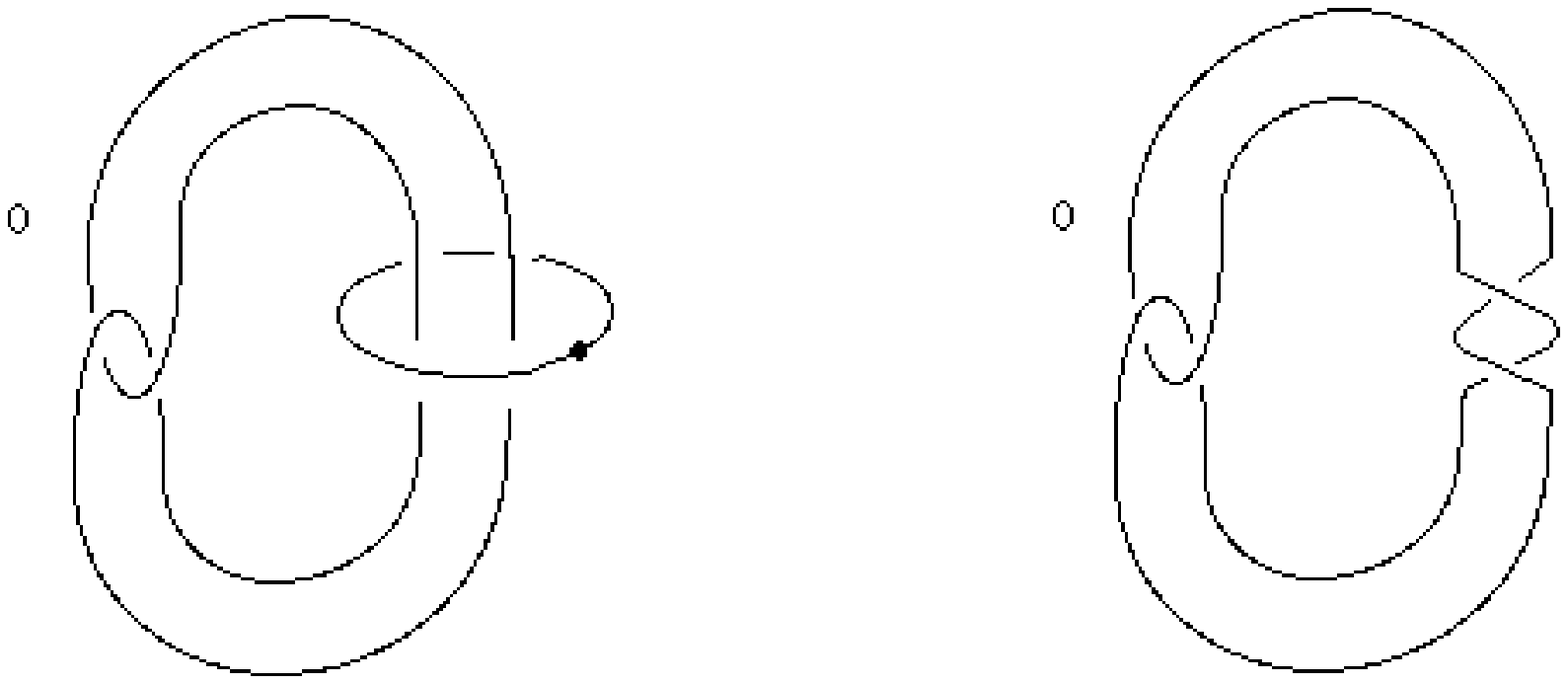}
\caption{}
\end{figure}

\begin{thm} There are compact smooth manifolds $F^{*} $ and
$Q^{*}$, which are homeomorphic but not diffeomorphic to $F$ and $Q$,
respectively. Also, $ F^{*}$ is obtained by removing a tubular
neighborhood of a properly imbedded $2$-disc 
$f: B^{2}\hookrightarrow S^{2}\times B^{2} $ from $S^{2}\times B^{2} $.
\end{thm}

In Figure 1 the circle $f(\partial B^{2})$ corresponds to the
 circle with dot. Existence of a fake cusp was
first established jointly with R.Matveyev as application of methods of
\cite{AM}

\section{Seiberg-Witten invariants} 

Let $X$ be a closed smooth $4$-manifold and $ Spin_{c}(X) $ be the set of
$ Spin_{c} $ structures on $ X $. 
In case $H_{1}(X)$ has no $2$-torsion $ Spin_{c}(X) $ can be identified by
 $$ Spin_{c}(X) =\{ a\in H^{2}(X,{\bf
Z})\;|\; a=w_{2}(TX) \; (\mbox{mod } 2) \}$$
Recall Seiberg-Witten invariant
$$SW_{X}:Spin_{c}(X) \to {\bf Z} $$
A classes $a\in H^{2}(X,{\bf Z})$ is called {\it basic} if $SW_{X}(a)\neq 0$.
It is known that there are finitely many basic classes and if $a$ is basic then
so is $-a$, and $$SW_{X}(-a)= (-1)^{\e}SW_{X}(a)$$
where $ \e = (e(X)+\s(X))\;/4 $, and $e(X), \s(X)$ denote Euler characteristic
and signature. If ${\it B}=\{\pm \a_1, \pm \a_2,..,\pm a_n\}$ are
the basic classes, by denoting  $ a_{0}=SW_{X}(0), a_{j}=SW_{X}(\a_{j}) $
and $t_{j}=exp(\a_{j})$
 Seiberg-Witten invariants can be assembled a single polynomial \cite{FS}
$$ SW_{X}= a_0 + \sum_{j=0}^{n}a_{i}(t_{j} + (-1)^{\e}t_{j}^{-1})$$

In \cite{FS} Fintushel and Stern introduced a method of modifying a
$4$-manifold by using a knot in $ S^{3}$, which changes its
Seiberg-Witten invariants without changing its homeomorphism type.
Let  $X$ be a closed smooth $4$-manifold and $ T^2\subset X $ be an
imbedded $2$-torus with trivial normal bundle. Assume that this torus
lies in in a cusp neighborhood; this means that inside of $X$ the
tubular neighborhood $T^2\times B^2 $ of the $2$-torus (first picture of
Figure 2) is contained in a cusp C (the second picture Figure 2). Call
such a $2$-torus in $X$ a
\textit{c-imbedded torus}

\begin{figure}[htb]
\includegraphics[scale=.5]{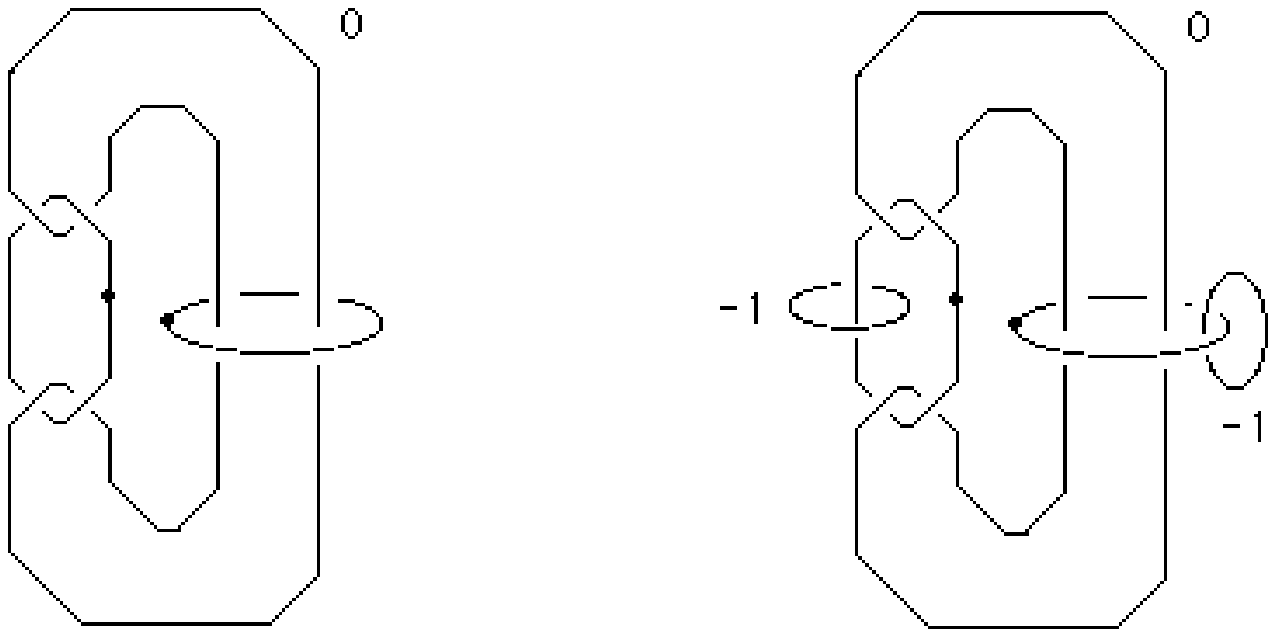}
\caption{}
\end{figure}

Let $ K \subset S^3 $ be a knot, and $ N=N(K)\approx  K\times B^2 $ be
the trivialization of its open tubular neighborhood given by the
$0$-framing. Let 
$ \varphi: \partial (T^{2}\times B^{2})\to \partial  (K\times B^{2}) \times
S^{1}
$ be any diffeomorphism with $ \varphi (p\times \partial B^2)=K\times p $
where
$ p\in T^{2} $ is a point. Define: 

$$X_{K}= (X-T^2\times B^2)
\smile_{\varphi} (S^{3}-N)\times S^1$$ 
 Let $[T]$ be the homology class in $H_{2}(X_{K};{\bf Z})$
induced from  $T^{2}\subset X$, and $ t=\mbox{exp}(2[T]) $, and let
$\Delta_{K}(t)$ be the Alexander polynomial of the knot $K$ (as a symmetric
Laurent polynomial), then

\begin{thm}{\cite{FS}}: Let $X$ be a smooth manifold as above, and 
$ K\subset S^{3} $ be a knot , then
Seiberg-Witten invariants of $ X_{K} $ can be computed
$$ SW_{X_{K}}=SW_{X} . \;\Delta_{K}(t) $$
\end{thm}

\section{Handlebody of $X_{K}$ }
Here we will give a general algorithm of describing 
the handlebody description of $X_{K}$ from the handles of $X$. It could
be beneficial for the reader to compare the  steps of this section with
\cite{A3}. Let
$ K\subset S^{3}
$ be a knot, depicted as the left handed trefoil knot in Figure 3, and $N$ be
its open tubular neighborhood. We claim that the linking circles
$\alpha $ and
$ \beta $ are the {\it core circles} of the $1$- handles of the
handlebody  of $ S^3-N $ (Heegard diagram).

\begin{figure}[htb]
\includegraphics[scale=.5]{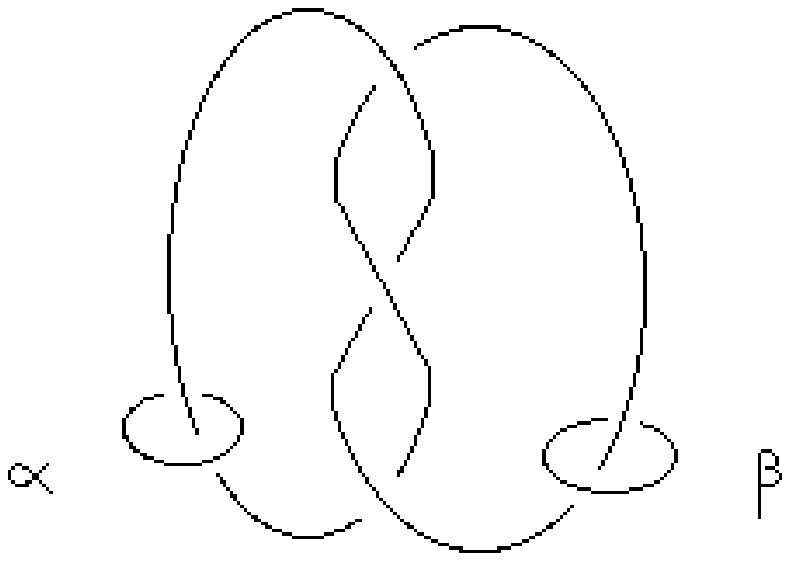}
\caption{}
\end{figure}

This can easily
checked by the process described in \cite{R} pp.250, which turns a
surgery description of a $3$-manifold to its Heegard diagram.
Steps in Figures 4 describes this process (i.e. attach canceling pair of
$1$ and $2$ - handles to $ S^{3}-N $ until the complement becomes a solid
handlebody). Write:

\begin{figure}[htb]
\includegraphics[scale=.5]{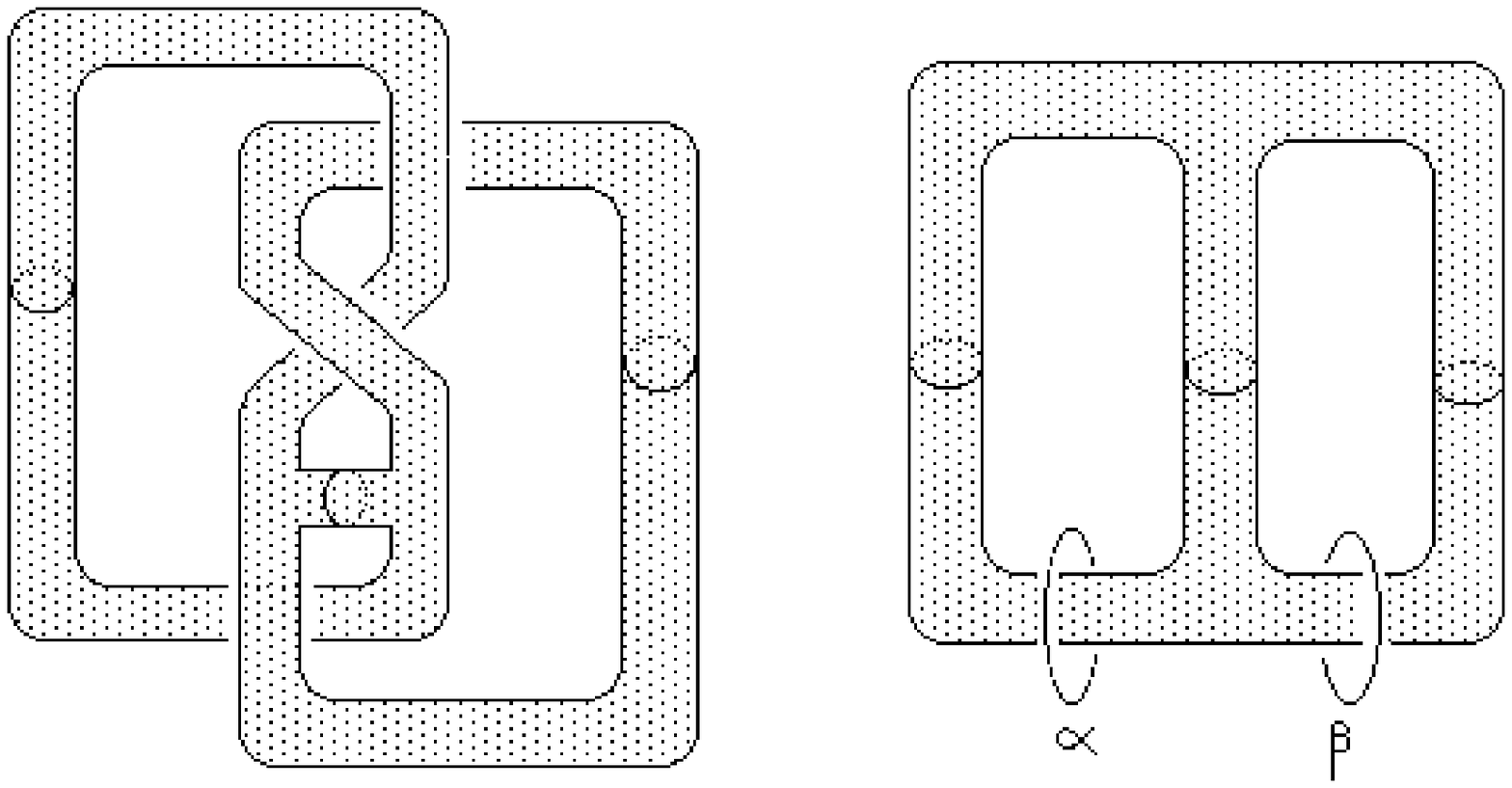}
\caption{}
\end{figure}

$$(S^3 -N) \times S^1= (S^3 -N)\times I_{+}
\smile (S^3 -N)\times I_{-}$$ 
where $I_{\pm} \approx I =[0,1]$
are closed intervals and the union is taken along the boundaries,
i.e. along $ (S^3 -N)\sqcup (S^3 -N)$. 
 Up to attaching a $3$-handle, $ (S^3 -N)\times I_{-}$ is
obtained  by removing the tubular neighborhood of a properly
imbedded arc (with trefoil knot tied on it) from $B^3$, and
crossing it with $I$ as indicated in Figure 5.

\begin{figure}[htb]
\includegraphics[scale=.6]{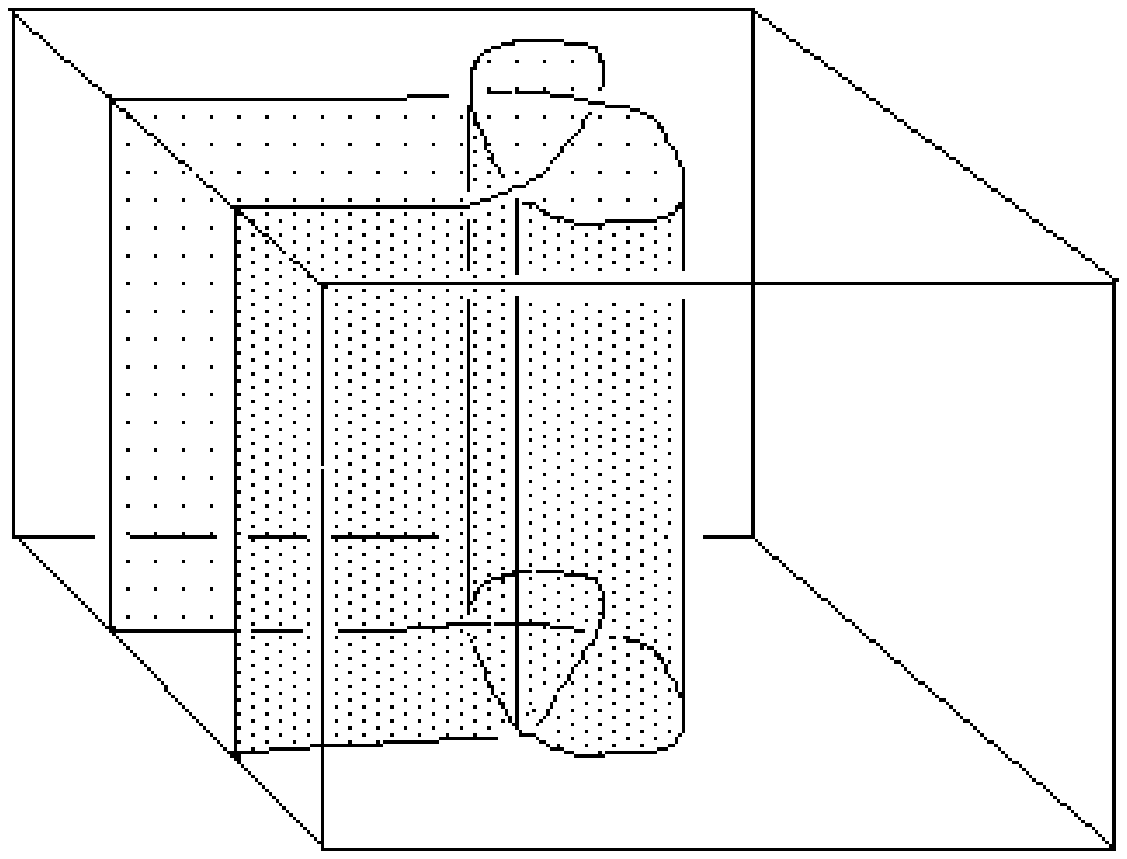}
\caption{}
\end{figure}

Equivalently, $ (S^3 -N)\times I_{-}$ is obtained by removing the ``usual"
slice disc from
$B^4$, which trefoil knot connected summed with its mirror image $K\# (-K)$
bounds, as indicated in Figure 6. The dot on the knot $\;K\#(-K)\;$
in Figure 6 indicates that the tubular neighborhood of the obvious slice
disc which it bounds is removed from $ B^4 $. We will refer this
as \textit{slice 1-handle}. This notation was discussed in
\cite{AK} and \cite{A3}. 

\begin{figure}[htb]
\includegraphics[scale=.5]{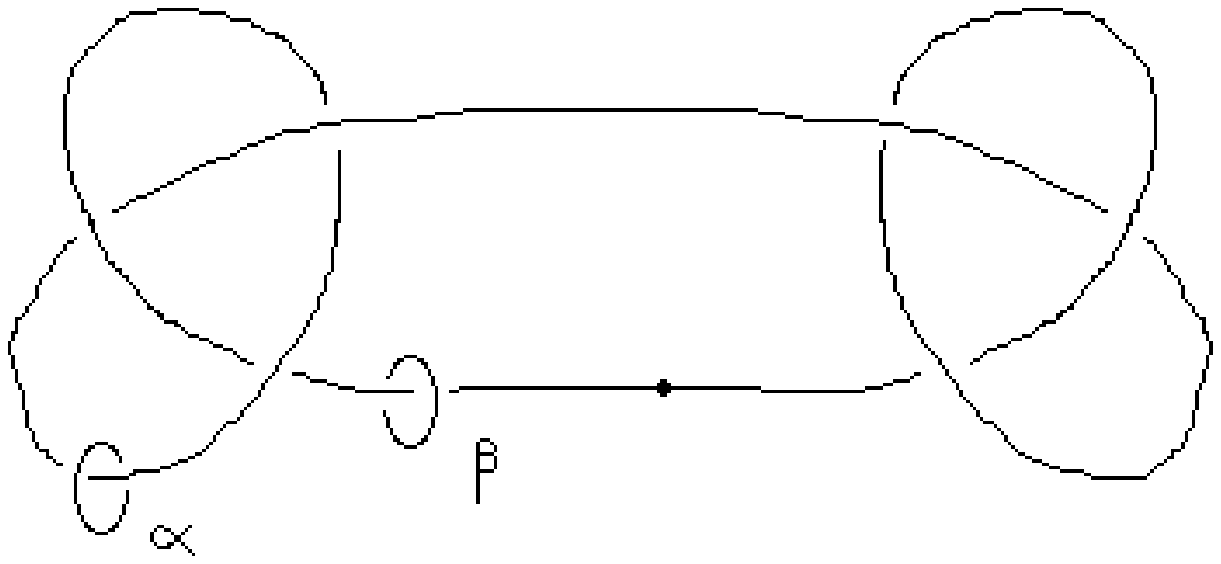}
\caption{}
\end{figure}

To get $ (S^3 -N)\times S^1 $ we glue the upside down
 handlebody of $ (S^3 -N)\times I_{+} $ to $ (S^3 -N)\times
I_{-} $. Clearly, up to attaching $3$-handles, this is achieved by
attaching $(S^3 -N)\times I_{-}$ one $1$-handle and two $2$-handles,
resulting from identification of the corresponding  $1$-handles $\a$
and
$\b$ (of the knot complements) in the two boundary components
of $(S^3-N)\times I_{-}$, as indicated in Figure 7.
So Figure 7 gives $( S^3-N)\times S^1$.

\begin{figure}[htb]
\includegraphics[scale=.5]{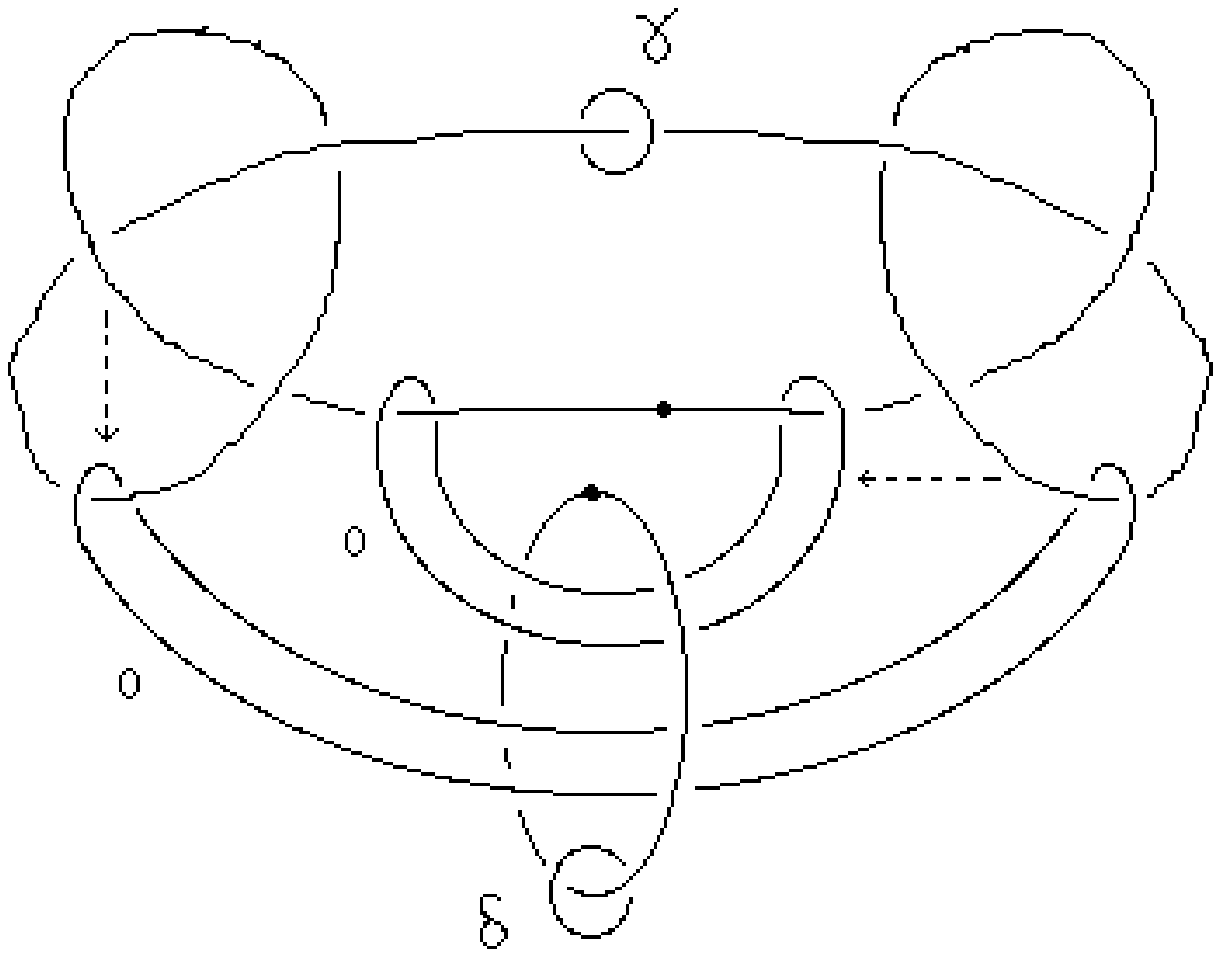}
\caption{}
\end{figure}

At this stage it is
instructive to check that the boundary of the handlebody of Figure
7 is indeed $T^3$. This can be seen by removing the ``dot" from
the ``slice 1-handle" then performing handle moves as indicated by
the arrows in Figure 7. This gives the first picture of Figure 8, then by
sliding one of the $0$-framed handles from the other one we get the
second picture Figure 8. Since the trivial knot with $0$-framing
is canceled by a $3$-handle, we see that Figure 8 is $ T^2\times B^2 $
(with boundary $T^3 $). Hence the reverse operation 
Figure 8 $\to$ Figure 7 corresponds the modification $X\to X_{K}$

\begin{figure}[htb]
\includegraphics[scale=.5]{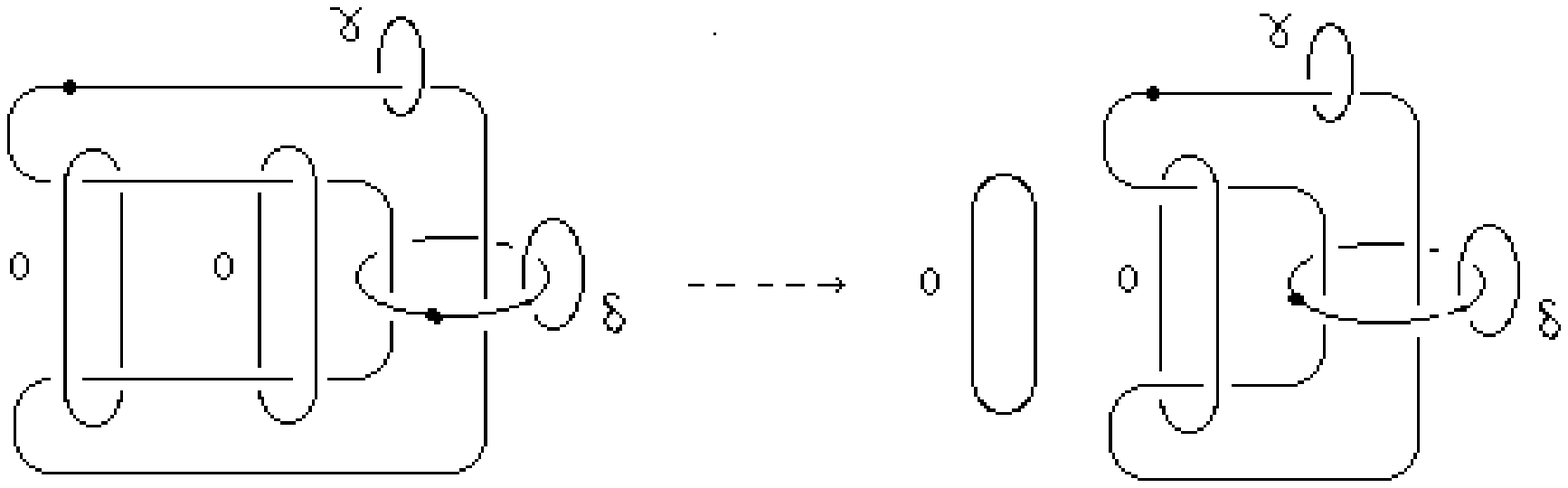}
\caption{}
\end{figure}

We can see the operation  Figure 8 $\to $ Figure 7 directly as
follows: Start with $ T^2\times B^2 $, by attaching a canceling
pair of $2$ and $3$ handles (i.e. by introducing an unknot with
$0$-framing), and then by sliding this new $2$-handle over the
$2$-handle of $T^2\times B^2$ we get another handle description of
$T^2\times B^2$ in Figure 9. Now by removing the ``dot" on the
$1$-handle (i.e. by surgery turning $1$-handle $S^1\times B^3$ to
$B^2\times S^2$) and performing the handle slides as indicated by the
arrows of Figure 9 (and putting a ``dot" on a resulting $0$-framed knot), we
get Figure 7. The ``circled $1/2 $ notation" on one of the arrows of Figure
9 means that when doing handle slide put one half-twist on the band.

\begin{figure}[htb]
\includegraphics[scale=.5]{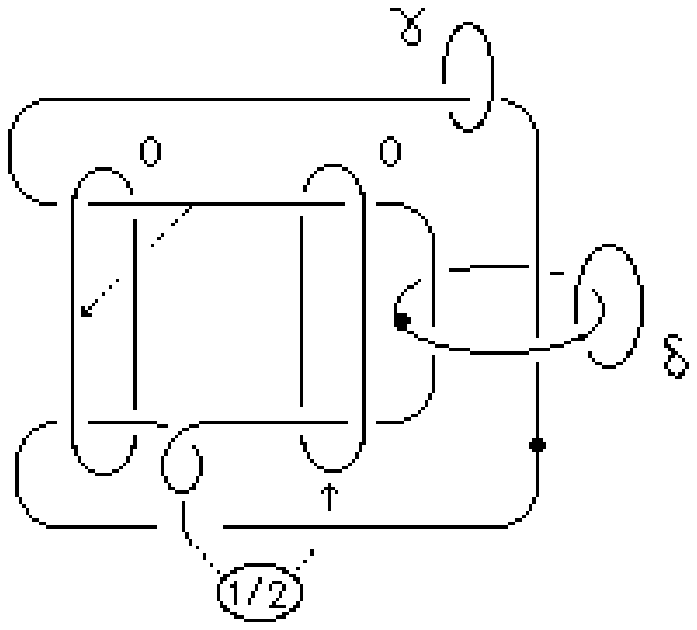}
\caption{}
\end{figure}

Hence to see exactly how the
operation $X\to X_{K}$ modifies the handles of $X$. We simply apply this
process to an imbedded $ T^2\times B^2 $ inside $X$, and trace the 
rest of the other handles of  $X$ along with it. For example, the operation 
Figure 9 $\to$ Figure 7 preserves the linking circles $\g$ and $\d$, as 
indicated in the figures. Therefore, if $T^2\times B^2$ lies in a
cusp neighborhood in $Q$ (i.e. if there are $-1$ framed 
handles attached to the knots $\g$ and $\d$ of Figure 9), Figure 9
becomes Figure 10, and the  operation $ X\to X_{K} $ corresponds to the
operation Figure 10 $\to$ Figure 11 (here disregard the loop $\t$ in 
Figures 10 and 11, it will be explained in the next paragraph).

\begin{figure}[htb]
\includegraphics[scale=.5]{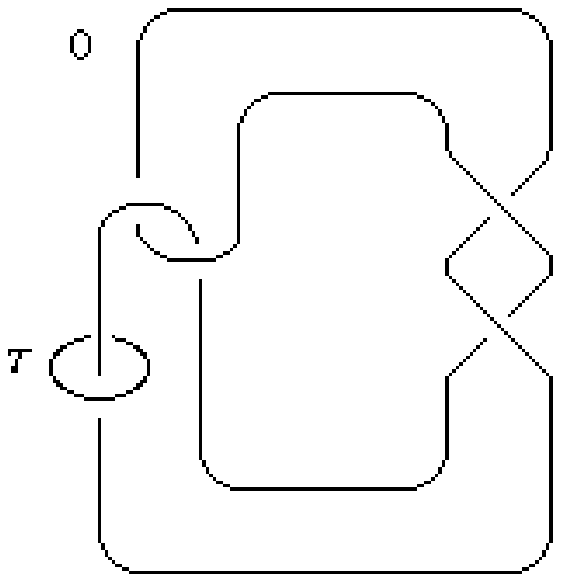}
\caption{}
\end{figure}

\begin{figure}[htb]
\includegraphics[scale=.5]{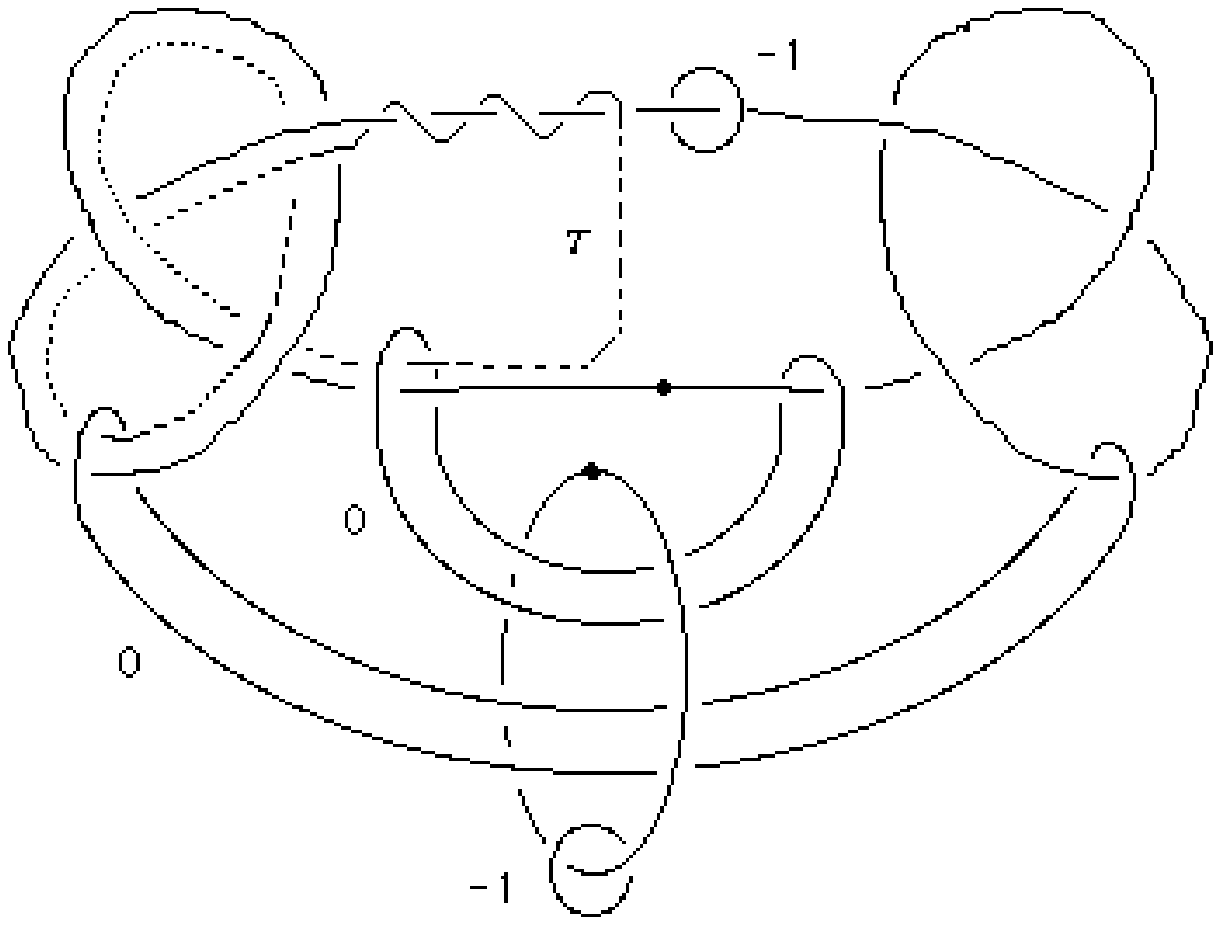}
\caption{}
\end{figure}

Warning: even though the operation Figure 9 $\to$ Figure 7 preserves
the loops $ \g $, $ \d $ of the $1$-handles, it does not preserve the
linking circle $\t$ of the $2$-handle of $ T^2\times B^2 $. Hence for
example, $\t$ of the cusp of Figure 10 is sent to the quite 
complicated loop of Figure 11 (also denoted by $\t$). 
Finally by drawing the slice $1$-handle of Figure 11 as two
$1$-handles and one $2$-handle , we see that
Figure 11 is diffeomorphic to Figure 12  (canceling one of the
$1$-handles of Figure 12 by the ``middle" $2$-handle gives the slice $1$-handle
of Figure 11, as in \cite{A3}).

\begin{figure}[htb]
\includegraphics[scale=.5]{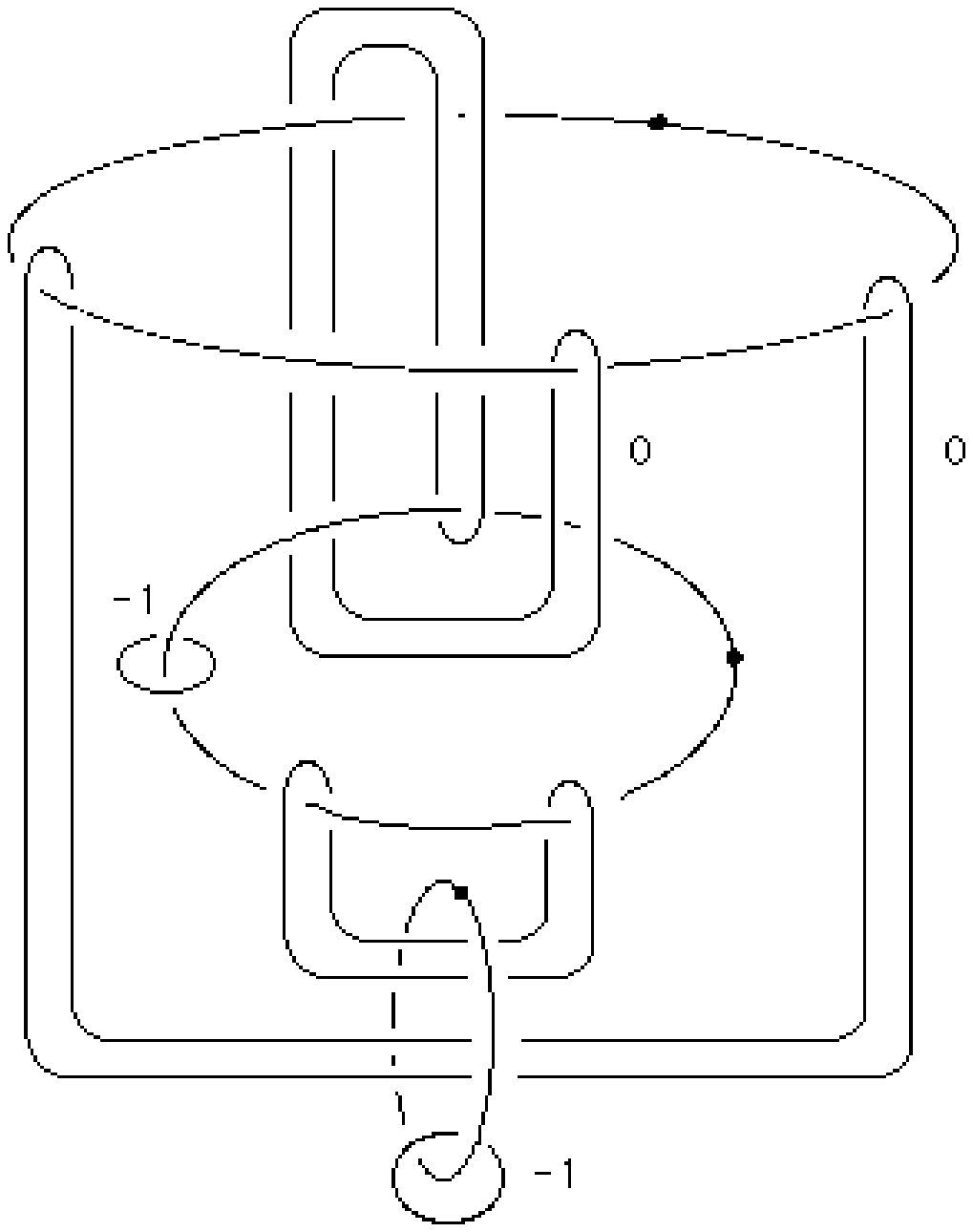}
\caption{}
\end{figure}

To sum up: The operation $X\to X_{K}$ changes an imbedded  $T^{2}\times
B^{2}$ inside of $ X $ by changing one of its $1$-handles with a  ``slice
$1$-handle" determined by the knot $K\# (-K) $, and
introducing 2-handles connecting the ``core circles" of $1$-handles of the
two knot complements (all except one) as in Figure 13.

\begin{figure}[htb]
\includegraphics[scale=.5]{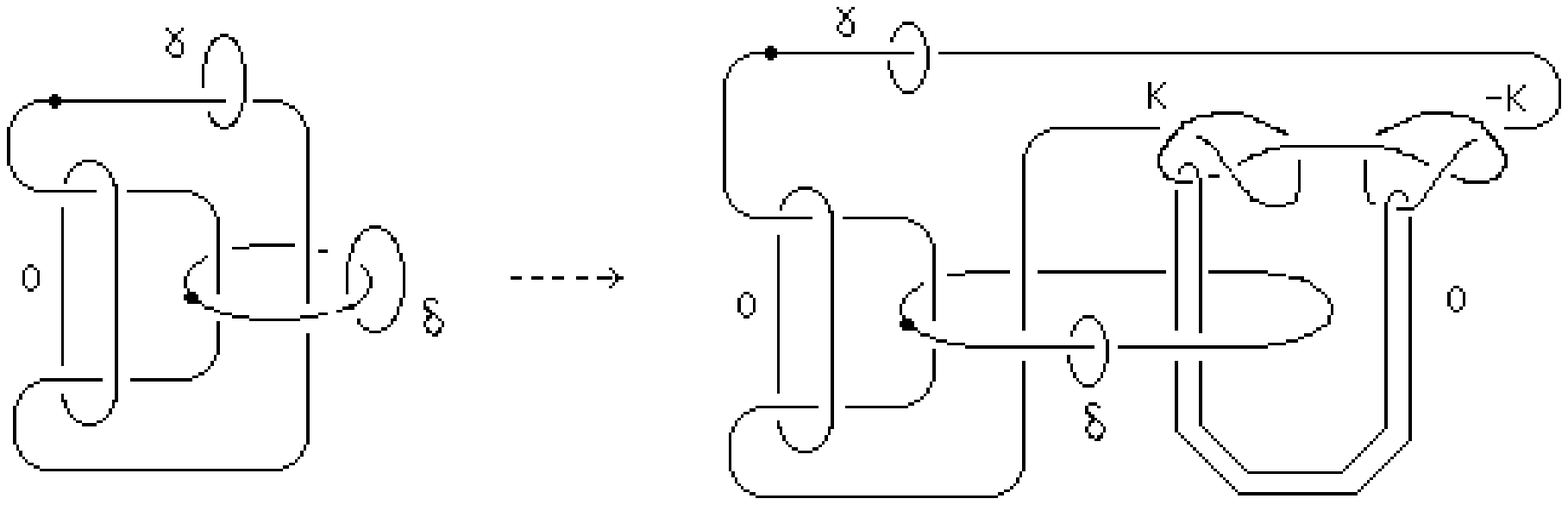}
\caption{}
\end{figure}

\section{Proof of the theorem }

Attaching $2$-handles with $-1$ framings to $\g $ and $\d $  of the first picture
of Figure 13  gives the cusp C. By Eliashberg's theorem $C$ is compact
Stein manifold (see \cite{AM} for brief review). C can be compactified
to a Kahler manifold $ X $; for example $ C $ sits in a K$3$ surface as a
codimension zero submanifold. By applying Theorem 1.1 to the torus 
$ T^{2}\subset C\subset X $, and a knot $ K $ is with nontrivial Alexander
polynomial, we obtain a fake copy  $ X_{K} $ of $ X $ (because $ X_{K} $
has different Seiberg-Witten invariant than $X$). Define $C^{*}=
C_{K}
\subset X_{K} $, then
$C^{*}$ can not be diffeomorphic to $ C $, otherwise the identity map 
$ id: X - \mbox{int}(C)\to X_{K}-\mbox{int}(C^{*}) $ would extend to a
diffeomorphism $ X\to X_{K} $. Recall that all self-diffeomorphisms of the
boundary of the cusp $ C $ extends to $ C $. Figure 12 is the handlebody of
$C^{*}$ (in case $K$ is the trefoil knot)

\begin{figure}[htb]
\includegraphics[scale=.5]{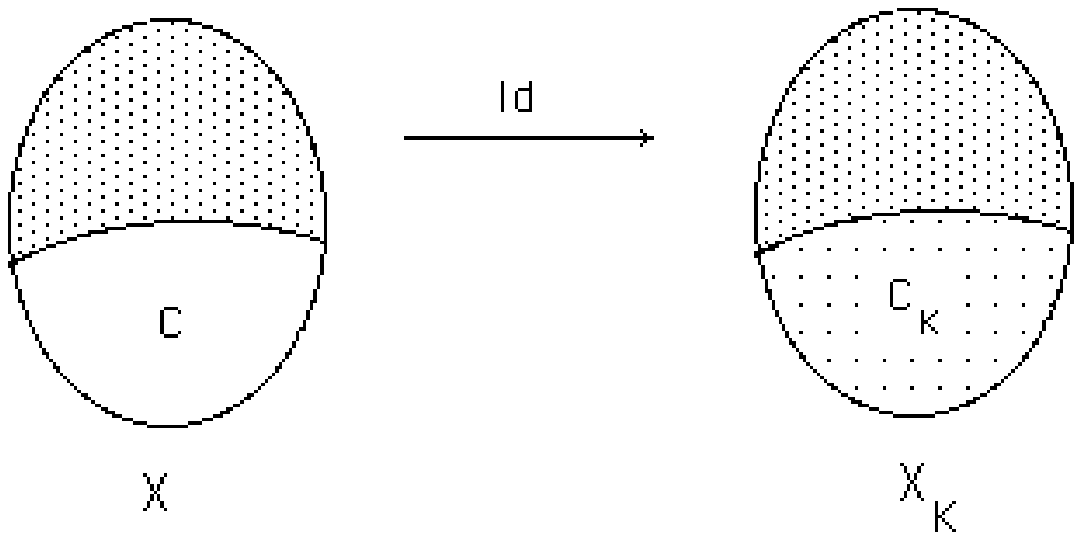}
\caption{}
\end{figure}

Attaching one $2$-handle with $-1$ framings to either one of the circles $\g $ or
$\d $  (say to $ \g $) of the first picture of Figure 13  gives the fishtail
$F$. We can think of of $ F $ being obtained from $ C $ by ``undoing" one of
its
$2$-handles, i.e. removing a thickened disc $ D $ from $C$ (dual $2$-handle of
$ \d $). The boundary 
$ \partial D $ corresponds the small trivially linking circle of the $-1$ framed
circle $ \d $ ; so removing $ D $ from $ C $ corresponds to putting a ``dot"
on this dual circle (and hence canceling the
$2$-handle $\d$ from $C$). By extending $ id: X - \mbox{int}(C)\to
X_{K}-\mbox{int}(C^{*}) $ across $D$ gives a diffeomorphism 
$ f: \partial F\to
\partial F_{K} $ which does not extend over the interior (otherwise $ X $
and $ X_{K} $ would be diffeomorphic).

\begin{figure}[htb]
\includegraphics[scale=.5]{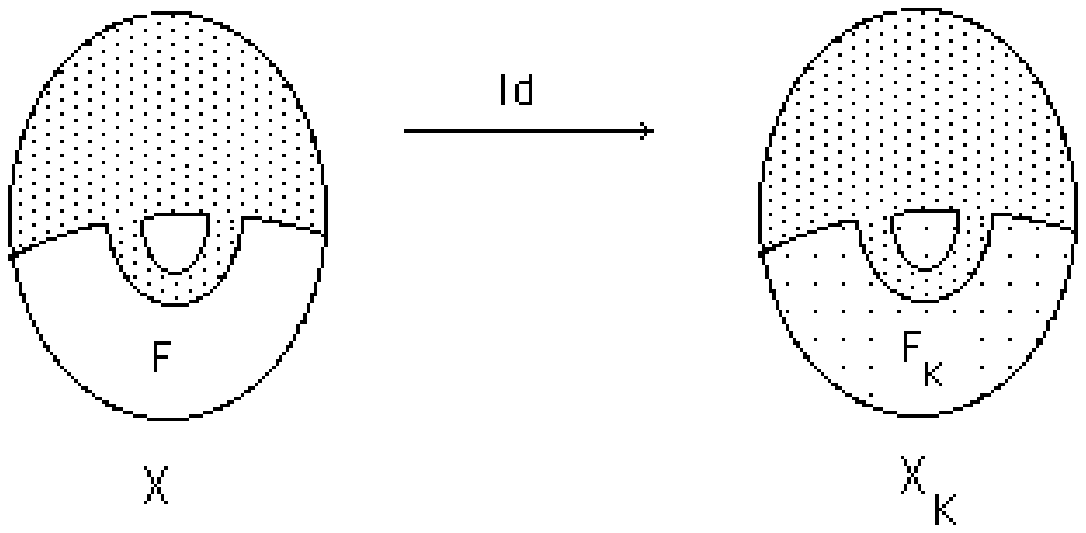}
\caption{}
\end{figure}

 Define $ F^{*} = F_{K} $, so removing
thickened $D$ from $C^{*}$ (Figure 11) gives $F^{*}$ (Figure 16). It is
easy to verify that $F^{*}$ is homotopy equivalent to $F$.  
Furthermore we can
verify that $F^{*}$ is obtained from $S^{2}\times B^{2}$ by removing an
imbedded disc
$D$ as follows: By the handle moves of Figure 7, we see that on the
boundary the position of the uknotted ``circle with dot" of Figure 16 is
the same as the ``circle with dot" of Figure 1 (i.e. F) ; also removing the
`circle with dot" from Figure 16 results $ S^{2}\times B^{2} $ (this can be
verified by going from Figure 16 to the handle presentation of Figure 12).

\begin{figure}[htb]
\includegraphics[scale=.5]{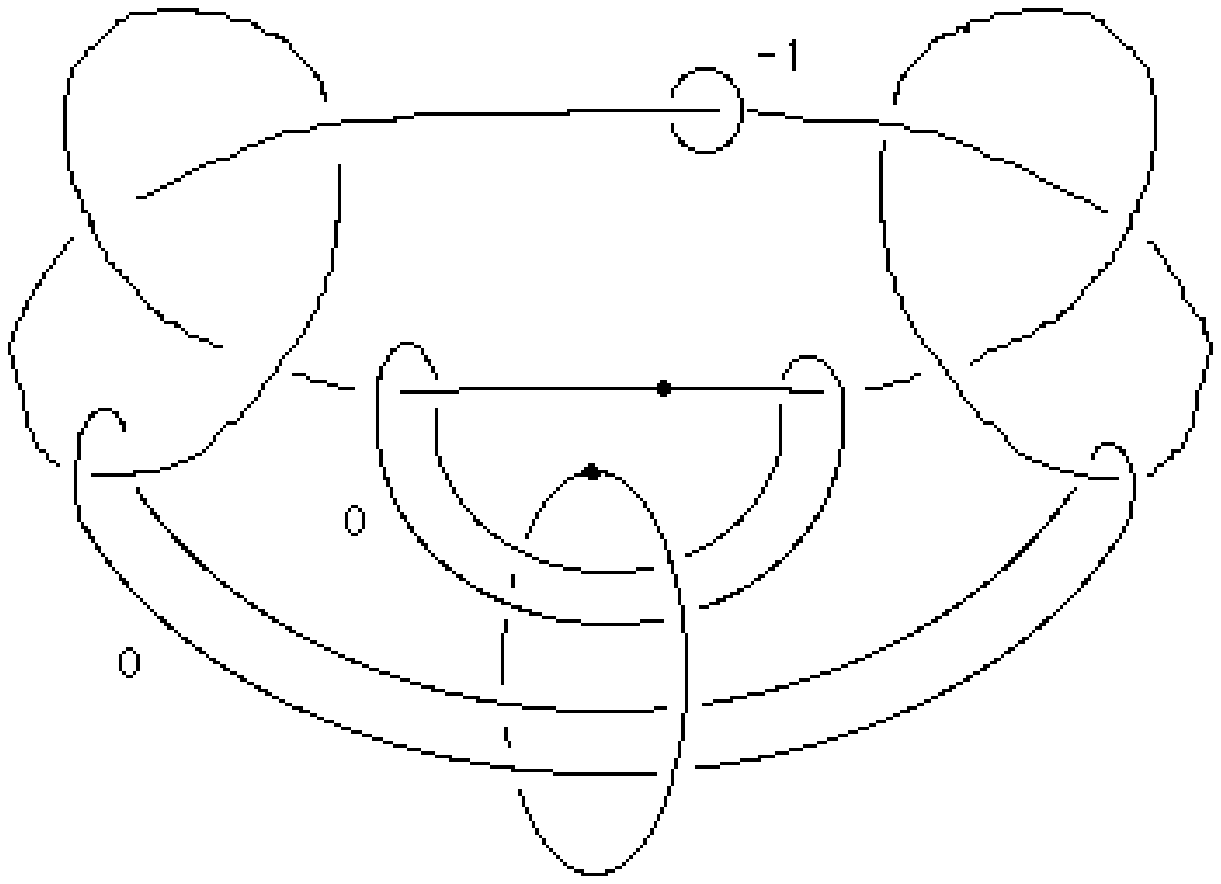}
\caption{}
\end{figure}

It remains to be check that $F^{*}$ is not diffeomorphic
to $F$; so far we only know a particular diffeomorphism of $f:\partial F\to
\partial F^{*}$ does not extend inside. Unfortunately unlike C, not every self
diffeomorphism of $ \partial F $ extends over $F$. But from the construction
it is easy to see that $f$ extends to a homotopy equivalence $F^{*}\to F$,
hence the following lemma implies that $F^{*}$ can not be diffeomorphic to
$F$.

\begin{lem} If a diffeomorphism $ f:\partial F\to \partial F $ extends to a
self homotopy equivalence it extends to a diffeomorphism $ F\to F $ 
\end{lem} 

\begin{proof} (outlined by R.Gompf): It is known that  $\partial F$ is
a
$T^{2}$ bundle over
$S^{1}$ with monodromy $ \left( \begin{array}{cc} 1&1\\ 0&1 \end{array}\right) $
\cite{K}.
 By standard $3$-manifold theory $f$ can be isotoped
to a fiber preserving isotopy. By composing with obvious diffeomorphism that
extends, we can assume that the fiber orientation is preserved. Since $f$ has to
commute with monodromy it fixes the ``vanishing cycle" C, corresponding to $
(1,0) $. So $f$ is a composition of Dehn twist along the horizontal torus $C\times
S^{1}$ along $S^{1}$ direction, and Dehn twist along the fiber in $C$ direction
(Dehn twist orthogonal to $C$ is ruled out since it does not extend to homotopy
equivalence $F\to F$),  all these diffeomorphisms extend to $F$. \end{proof}.

\section{Exotic knottings of the cusp and the fishtail}

We can describe $S^4$ as a union of two fishtails along the boundary, and
$S^{2}\times S^{2}$ as a union of two cusps along the boundary. 
 Figure 17 describe these
identifications. 
For example, attaching an upside-down copy of $-F$ to $F$ has
an affect of attaching a $2$ and $3$-handles to $F$ as described in the
first part of Figure 17: Attaching the $2$-handle to $F$ gives
$S^{2}\times B^{2}$ which, after attaching a $3$-handle, becomes $S^{4}$.
Hence we have imbeddings of singular $2$-spheres $f_{0}:F\to S^{4}$ and $
g_{0}:C\to S^{2}\times S^{2} $

\begin{figure}[htb]
\includegraphics[scale=.6]{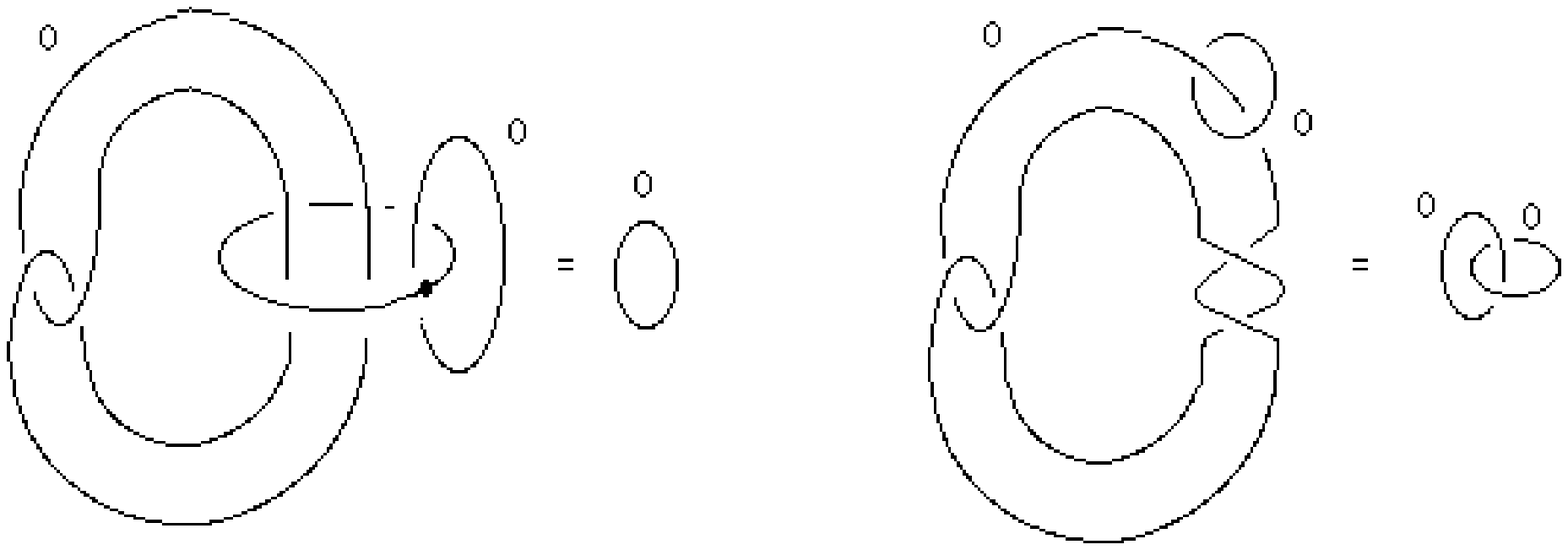}
\caption{}
\end{figure}

\begin{thm} There are imbeddings 
$ f_{1}:F\to S^{4} $ and $ g_{1}:C\to S^{2}\times
S^{2} $ that are topologically isotopic but not smoothly isotopic  to the
imbeddings $f_{0}$ and $g_{0}$ 
\end{thm}

\begin{proof} By replacing one of the fishtails in
$S^{4}=F\smile_{\partial} -F$ by $F^{*}$ we obtain a homotopy $4$-sphere, which
can be easily checked to be $S^{4} = F^{*}\smile_{\partial} -F $ (Figure 18).
Similarly by replacing one of the cusps in $ S^{2}\times
S^{2}= C\smile_{\partial}-C $ we obtain a homotopy $S^{2}\times S^{2}$ which can be
checked to be the standard $ S^{2}\times S^{2}= C^{*}\smile_{\partial}-C $.
Unlike the previous case, this check is surprisingly difficult (it requires the
proof of Scharlemann's conjecture). 

Figure 19 is the handlebody of
$ C^{*}\smile_{\partial}-C $, it is diffeomorphic to Figure 20 (to see
this, in Figure 20 slide one of the small $-1$ circles over one of the $0$-framed
handles, going through the $1$-handle, then slide $-1$ framed $2$-handle over it).
To identify Figure 20 by $ S^{2}\times S^{2}$ we need to first recall the
handlebody picture of $ \Sigma \times S^{1} $, where $\Sigma $ is the Poincare
homology sphere (\cite{A3}): Figure 21. From
\cite{A3} we know that surgering $ \Sigma \times S^{1} $ (along the loop trivially
linking the slice 1-handle) gives $ S^{3}\times S^{1}\# S^{2}\times S^{2} $, and
surgering once more the obvious $S^{1}$ gives $S^{2}\times S^{2}$. Performing
these two surgeries corresponds to introducing pair of $0$ and $-1$ framed two
handles as indicated in Figure 22. By canceling two
$2$ and $3$-handle pairs from  Figure 22 gives Figure 20. Via the handle moves of
Figure 7 one can check that, introducing the two little $0$-framed handles to
Figure 20, to obtain  Figure 22, has the affect of changing the boundary from $
S^{1}\times S^{2} $ to 
$ S^{1}\times S^{2}\# S^{1}\times S^{2} \# S^{1}\times S^{2} $ \end{proof}.

\begin{rem} The Fintushel-Stern operation $X\to X_{K}$ 
can be generalized by 
$$X\to X_{T}=(X-T^2\times B^2) \smile_{\varphi} (S^{3}\times S^1 - N_{T})$$
where  $N_{T}$ is an open tubular neighborhood of an imbedded $ T^{2} \subset
S^{3}\times S^{1}$. Surprisingly it turns out that this operation does not always
change the smooth structure of
$X$ (in particular it does not change the Seiberg-Witten invariant of $X$). 
Another generalization of this operation is by removing an open tubular
neighborhood of a Klein bottle $ N(F) $ (a twisted $ B^{2} $-bundle over $ F $)
from $ X $, and replacing it with a $ S^{3}-N(K) $ bundle over $ S^{1} $,
where $K\subset S^{3}$ is an invertible knot and $\psi : S^{3}- N(K)\to S^{3}-
N(K)$ is the inversion
$$X\to X_{K}=(X-N(F)) \smile_{\partial} (S^{3}- N(K))\times _{\psi} S^{1}$$
These operations will be studied in \cite{A4}.
\end{rem}

\begin{figure}[htb]
\includegraphics[scale=.44]{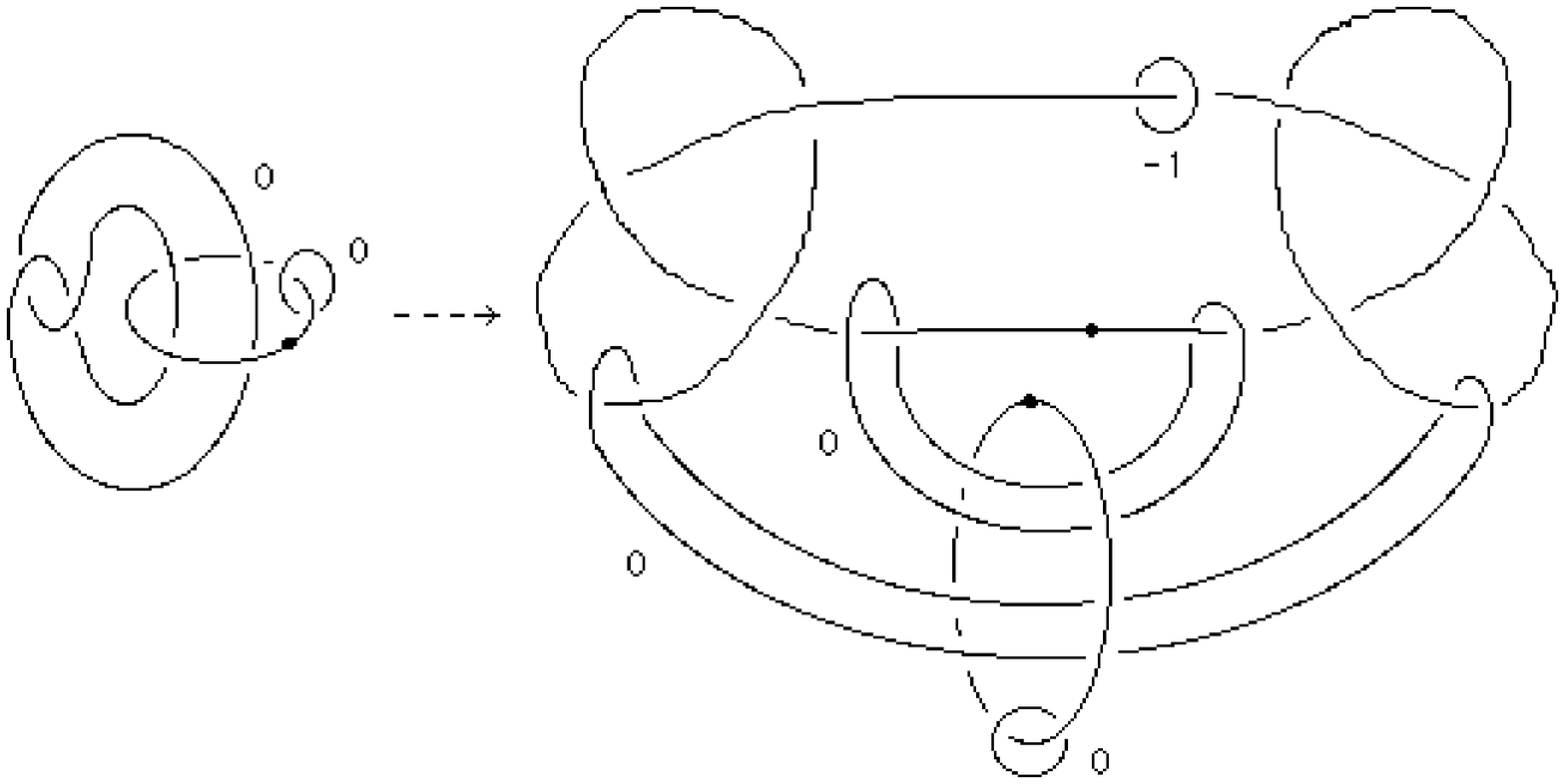}
\caption{}
\end{figure}

\begin{figure}[htb]
\includegraphics[scale=.45]{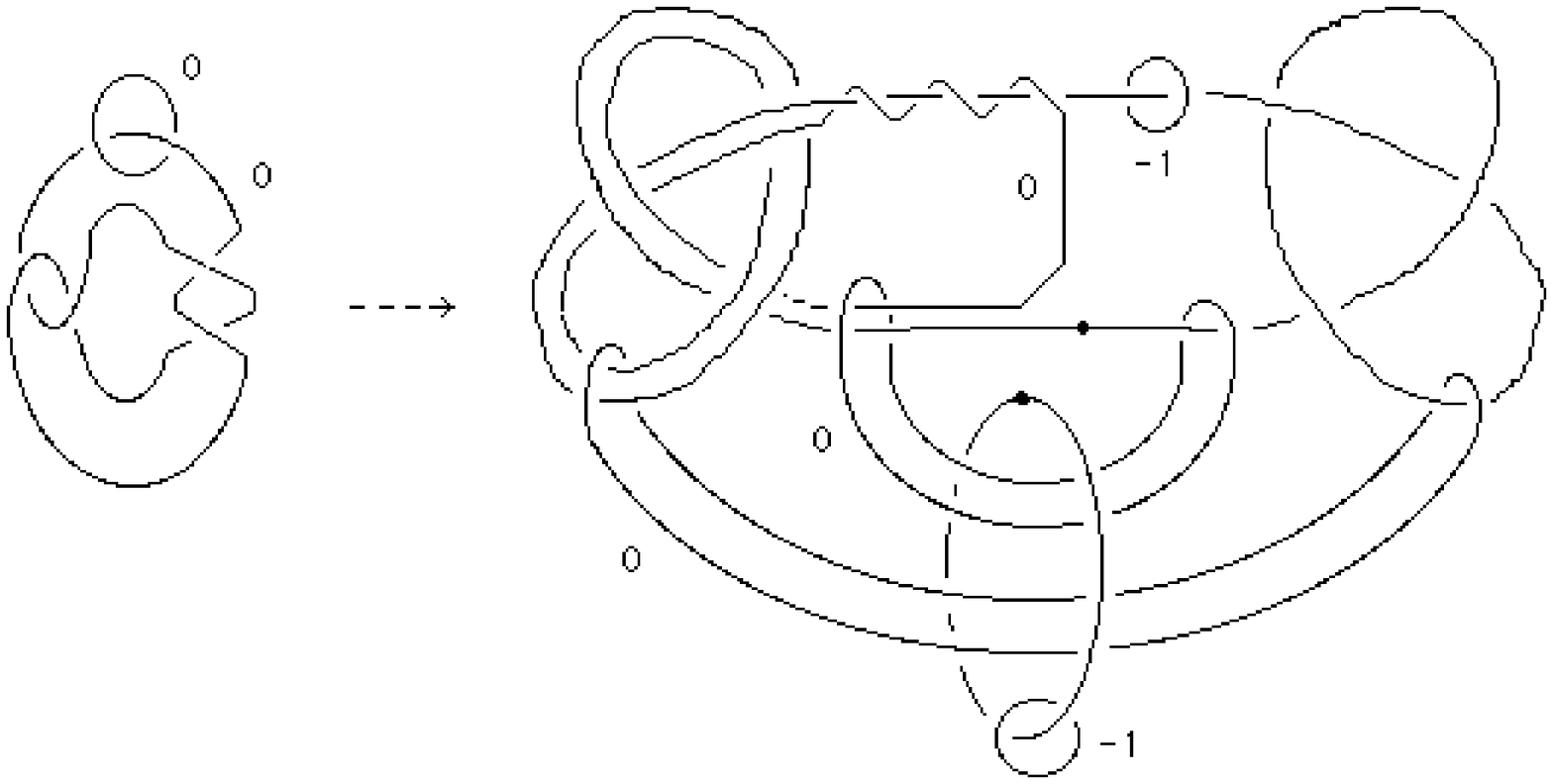}
\caption{}
\end{figure}

\begin{figure}[htb]
\includegraphics[scale=.45]{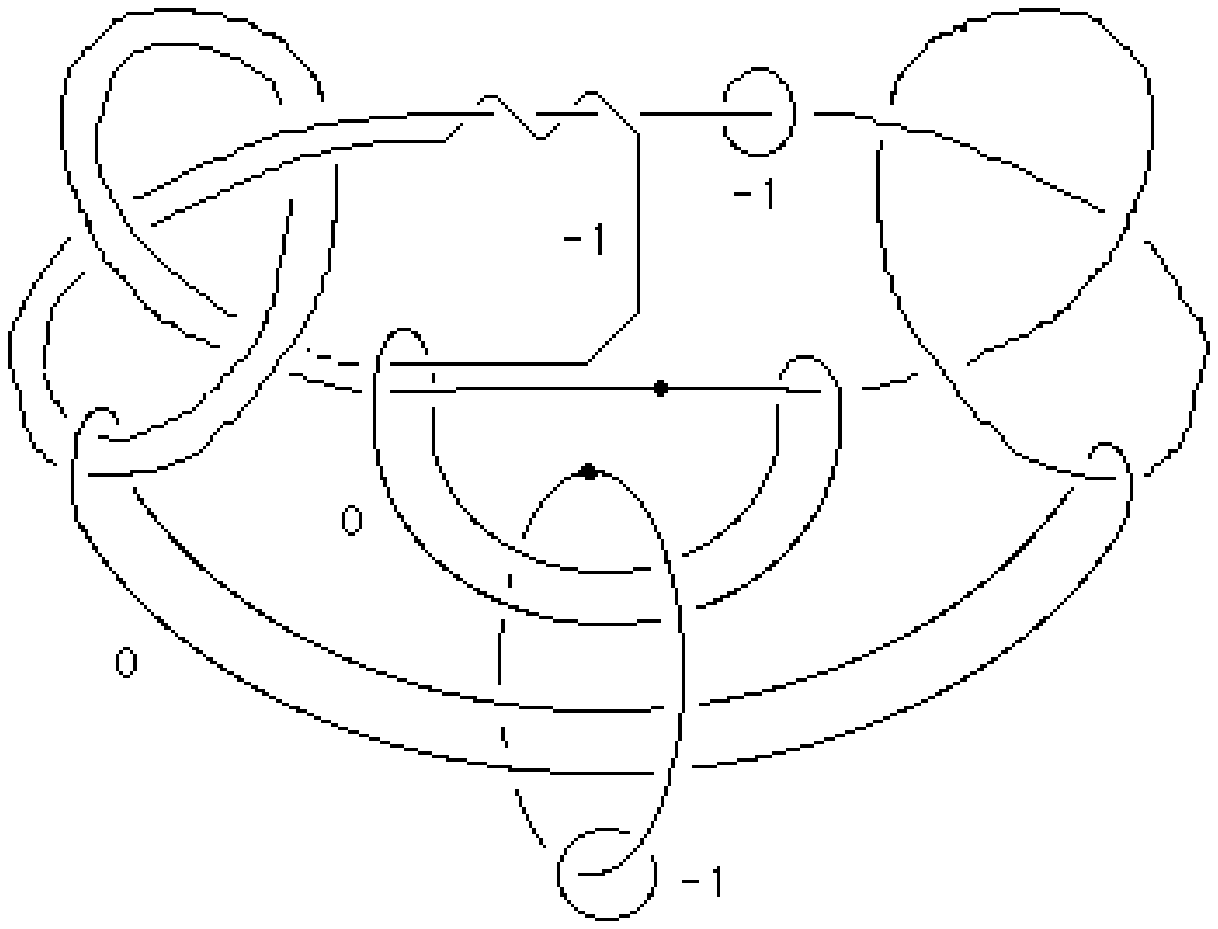}
\caption{}
\end{figure}

\begin{figure}[ht]
\includegraphics[scale=.5]{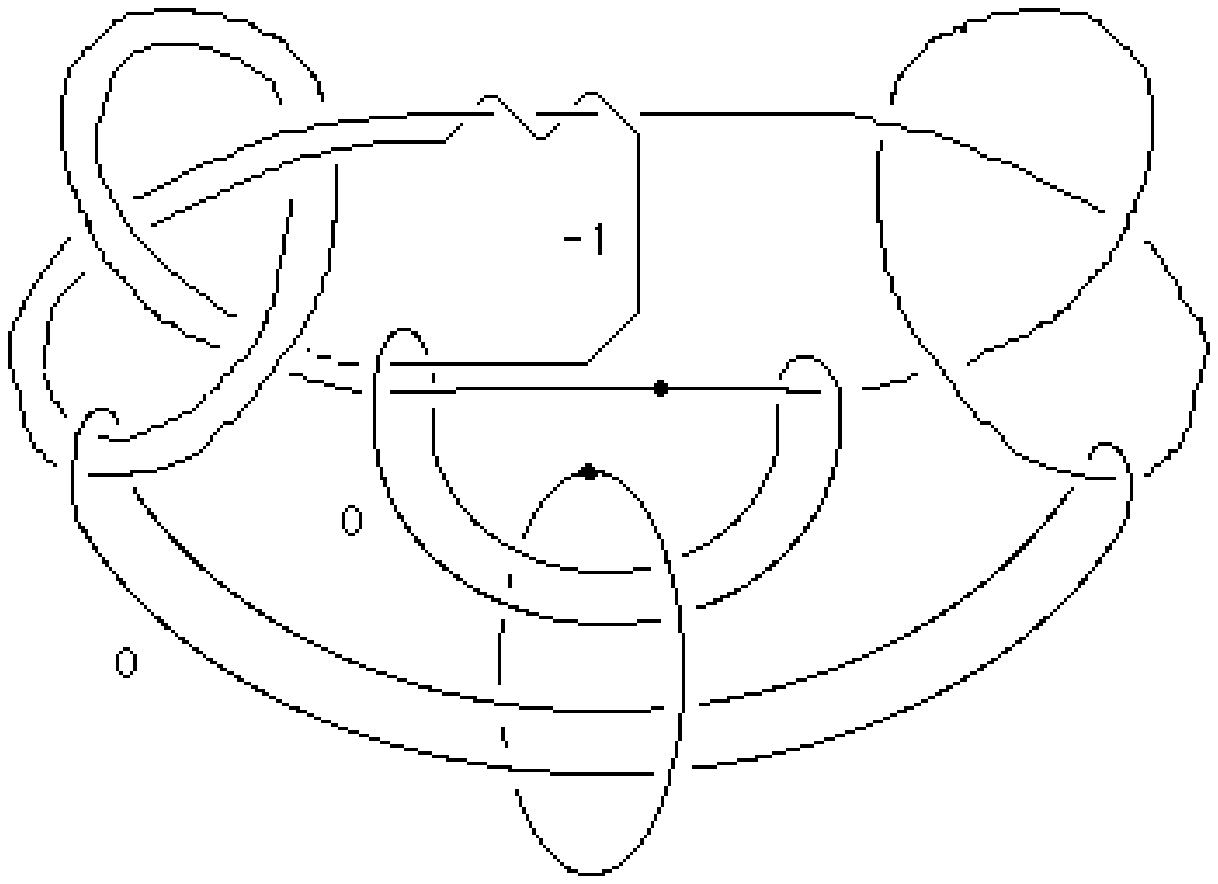}
\caption{}
\end{figure}

\begin{figure}[ht]
\includegraphics[scale=.5]{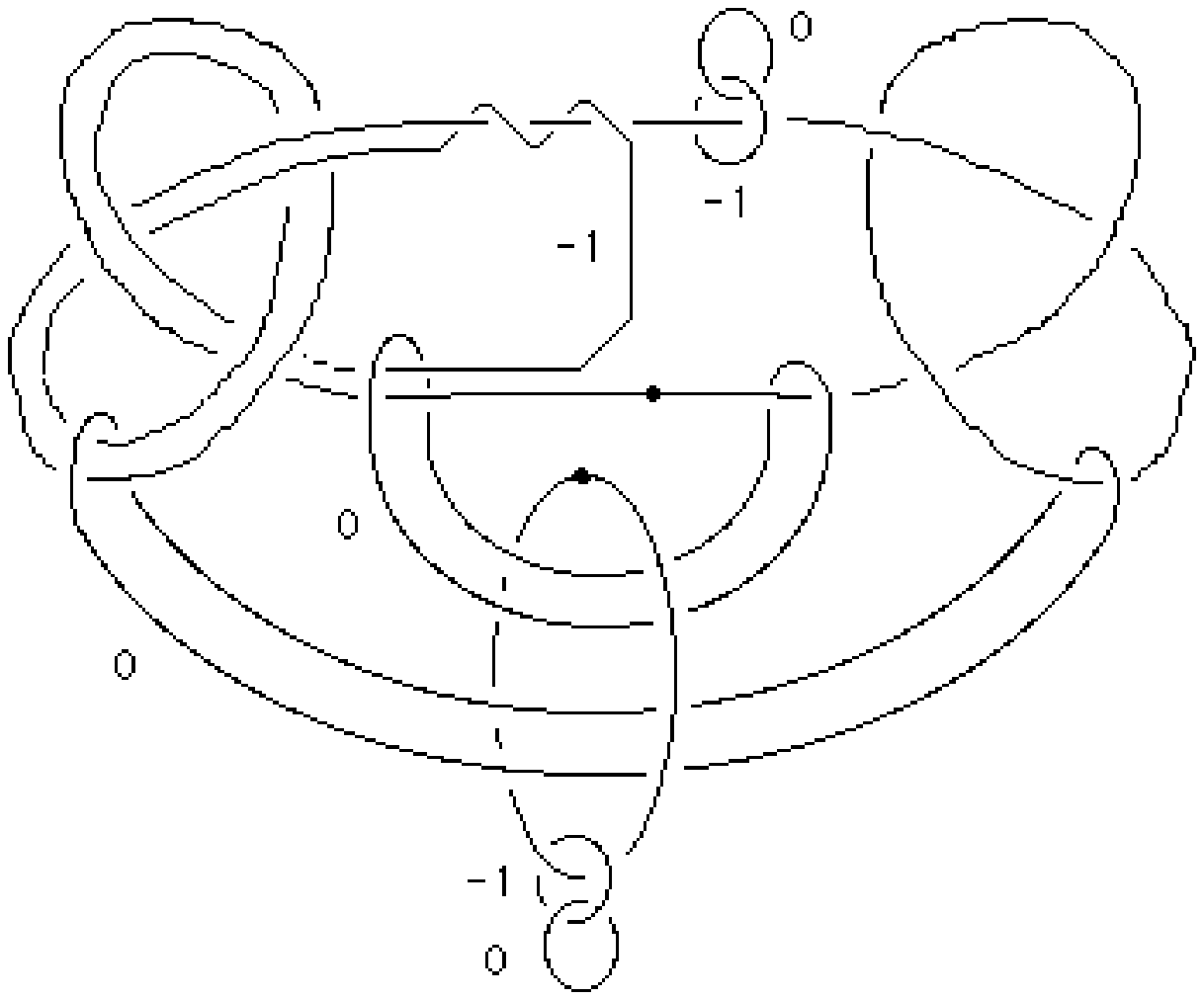}
\caption{}
\end{figure}

\end{document}